\documentclass[12pt]{amsart}
\usepackage{amsfonts,amssymb,amsthm,amsmath, color}
\usepackage{graphicx}
\usepackage{hyperref}

\newtheorem{thm}{Theorem}
\newcommand{\smat}[4]{\left(\begin{smallmatrix}
                 #1 & #2\\
                 #3 & #4
\end{smallmatrix}\right)}
\newcommand{\Z}{\mathbb{Z}}
\newcommand{\C}{\mathbb{C}}
\renewcommand{\H}{\mathbb{H}}
\newcommand{\abs}[1]{\left|#1\right|}
\newcommand{\SL}{{\text {\rm SL}}}

\numberwithin{equation}{section}

\begin{document}

\title{Zeros of weakly holomorphic modular forms of levels 2 and 3}
\author{Sharon Anne Garthwaite}
\address{Department of Mathematics, Bucknell University, Lewisburg,
PA 17870} \email{Sharon.Garthwaite@bucknell.edu}

\author{Paul Jenkins}
\address{Department of Mathematics, Brigham Young University, Provo,
UT 84602} \email{jenkins@math.byu.edu}

\begin{abstract}  Let $M_k^\sharp(N)$ be the space of weakly holomorphic
modular forms for $\Gamma_0(N)$ that are holomorphic at all cusps
except possibly at $\infty$.  We study a canonical basis for
$M_k^\sharp(2)$ and $M_k^\sharp(3)$ and prove that almost all
modular forms in this basis have the property that the majority of
their zeros in a fundamental domain lie on a lower boundary arc of
the fundamental domain.
\end{abstract}


\maketitle

\section{Introduction}

In studying a complex-valued function, it is natural to attempt to
locate the zeros of the function; in fact, one of the most famous
unsolved problems in mathematics asks whether the nontrivial zeros
of the Riemann zeta function $\zeta(s)$ lie on a particular line. In
this paper, we study the locations of the zeros of certain modular
forms, and show that most of their zeros in a fundamental domain
occur on a particular circular arc.

For the Eisenstein series, perhaps the easiest examples of modular
forms, a great deal is known about the locations of the zeros. In
the 1960s, Wohlfahrt~\cite{Wohl} showed that for even $4 \leq k \leq
26$, all zeros of the Eisenstein series $E_k(z)$ in the standard
fundamental domain for $\SL_2(\Z)$ lie on the unit circle $|z|=1$,
and R. A. Rankin~\cite{Ra3} extended the range of values of $k$ for
which this holds. Shortly afterward, F.K.C. Rankin and
Swinnerton-Dyer~\cite{RS-D} proved this result for all weights $k
\geq 4$.  R.A. Rankin~\cite{Ra1} also obtained the result for
certain Poincar\'{e} series, which generalize Eisenstein series.
Similar results have been obtained for Eisenstein series for
$\Gamma_0^*(2)$ and $\Gamma_0^*(3)$ by Miezaki, Nozaki, and
Shigezumi~\cite{MNS}, for Eisenstein series for $\Gamma_0^*(5)$ and
$\Gamma_0^*(7)$ and for Poincar\'{e} series for $\Gamma_0^*(2)$ and
$\Gamma_0^*(3)$ by Shigezumi~\cite{Shi2, Shi}, and for a family of
Eisenstein series for $\Gamma(2)$ by the first author, Long,
Swisher, and Treneer~\cite{GLST}.

The above results which locate the zeros of Eisenstein series and
Poincar\'e series use the same general idea of approximating the
modular form by an elementary function having the required number of
zeros on the arc. For example, the Eisenstein series $E_k(z)$ may be
written as a sum over an integer lattice.  When $z$ is restricted to
the unit circle, so that $z = e^{i\theta}$, the four main terms of
this series combine to give $2 e^{\frac{-ik\theta}{2}}
\cos\left(\frac{k\theta}{2}\right)$.  Rankin and Swinnerton-Dyer's
proof shows that the additional terms are small, so the zeros of
$E_k(z)$ are close to the zeros of this trigonometric function.

In 1997, Asai, Kaneko, and Ninomiya~\cite{AKN} used this idea to
study the zeros of polynomials related to the modular function
$j(\tau)$. The $j(\tau)$ function is a Hauptmodul, or an isomorphism
from the quotient of the upper half plane $\H$ under the action of
$\SL_2(\Z)$ to the complex plane $\C$.  It generates all modular
functions on $\SL_2(\Z)$, and it also parameterizes the isomorphism
classes of elliptic curves over $\C$. The image of $j(z)$ under the
Hecke operator $T_n$ is a polynomial in $j(z)$, which we write as
$P_n(j(z))$. Letting $\Delta(z)$ be the modular discriminant, which
is a weight 12 cusp form on $\SL_2(z)$ with no zeros in $\H$, a
generating function for the $P_n(j(z))$ is given by
\begin{equation}
\label{Jgen} \sum_{n=0}^\infty P_n(j(z))e^{2\pi i n \tau} =
\frac{E_4(\tau)^2E_6(\tau)}{\Delta(\tau)}\cdot\frac{1}{j(\tau)-j(z)}.
\end{equation}
Using this generating function, Asai et al.\ approximated the
polynomials $P_n(j(z))$ by trigonometric functions well enough to
prove that their zeros in the fundamental domain lie on the unit
circle.

Duke and the second author~\cite{DJ1} extended the results on $j(z)$
to a two-parameter family of modular forms that form bases for
spaces of weakly holomorphic modular forms of level 1.  In this
case, the connection to elementary functions is less direct.
Cauchy's integral formula relates the modular forms to a contour
integral of a generalized version of the generating function
(\ref{Jgen}).  An application of the residue theorem produces the
elementary functions, and the proof follows by bounding the integral
over a range of values for $\tau$ and $z$.  The zeros again lie on
the unit circle for many of the forms in the family, though in
contrast to previous results, it is known that this property does
not hold for all of the modular forms.

We mention one further result using a different technique by
Hahn~\cite{Ha}, who obtained general results on the zeros of
Eisenstein series for genus zero Fuchsian groups; the general idea
is an analogue of the classical argument that shows that the zeros
of an orthogonal polynomial all lie on the real line.

This question of locating zeros of modular forms is made even more
interesting by results of Rudnick~\cite{Ru} that showed that the
zeros of Hecke eigenforms of weight $k$, in a sense the orthogonal
complement of the Eisenstein series, are expected to become
equidistributed in the fundamental domain as $k \rightarrow \infty$; this conjecture was
proved by Holowinsky and Soundararajan~\cite{HS} as a consequence of
more general work on mass equidistribution for Hecke eigenforms.
Ghosh and Sarnak~\cite{GS} gave a lower bound for the density of
zeros lying on certain arcs for such eigenforms.  In a different
direction, Basraoui and Sebbar \cite{BS} proved that the
quasi-modular form $E_2(\tau)$ has infinitely many zeros that are
inequivalent under $\SL_2(\Z)$, and that none of these lie within
the fundamental domain.

In this paper, we examine modular forms in a basis for certain
spaces of weakly holomorphic modular forms of arbitrary integral
weight and levels $2$ and $3$. We show that for almost all of the
basis elements, most of their zeros in a fundamental domain for
$\Gamma_0(2)$ or $\Gamma_0(3)$ lie on a circular arc along the lower
boundary of the fundamental domain.  This is possible because we can
again approximate these modular forms by elementary functions;
however, the shape of the fundamental domain makes it difficult to
accurately locate all of the zeros.

\section{Definitions and statement of results} \label{defs}
Let $M_k(2)$ be the space of holomorphic modular forms of weight $k$
for the group $\Gamma_0(2) = \{\smat{a}{b}{c}{d} \in \SL_2(\Z) : c
\equiv 0 \pmod{2} \}$, and let $M_k^!(2)$ be the corresponding space
of weakly holomorphic modular forms, or modular forms that are
holomorphic on the upper half plane and meromorphic at the cusps.
Let $M_k^\sharp(2)$ be the subspace of $M_k^!(2)$ consisting of
forms which are holomorphic away from the cusp at $\infty$.  This
space appears, for instance, in~\cite{MP}, where it is shown that
traces of negative integral weight forms in such a space appear as
coefficients of certain half integral weight forms of level $4N$.
Modular forms in $M_k^\sharp(2)$ have been studied by
Ahlgren~\cite{Ahlgren}, who gave explicit formulas for the action of
the $\theta$-operator on forms in these spaces and obtained formulas
for the exponents of their infinite product expansions, and by
Andersen and the second author~\cite{JA}, who gave congruences for
the coefficients of a basis for $M_0^\sharp$.

For the group $\Gamma_0(2)$, we use a fundamental domain in the
upper half plane bounded by the lines $\textrm{Re}(z) =
-\frac{1}{2}$ and $\textrm{Re}(z) = \frac{1}{2}$ and the circles of
radius $\frac{1}{2}$ centered at $z = -\frac{1}{2}$ and $z =
\frac{1}{2}$.  We include the boundary on the left half of this
fundamental domain, which is equivalent to the opposite boundary
under the action of the matrices $\smat{1}{1}{0}{1}$ and
$\smat{1}{0}{2}{1}$.  The cusps of this fundamental domain may be
taken to be at $\infty$ and at $0$.

\begin{center}
\includegraphics[height=3in]{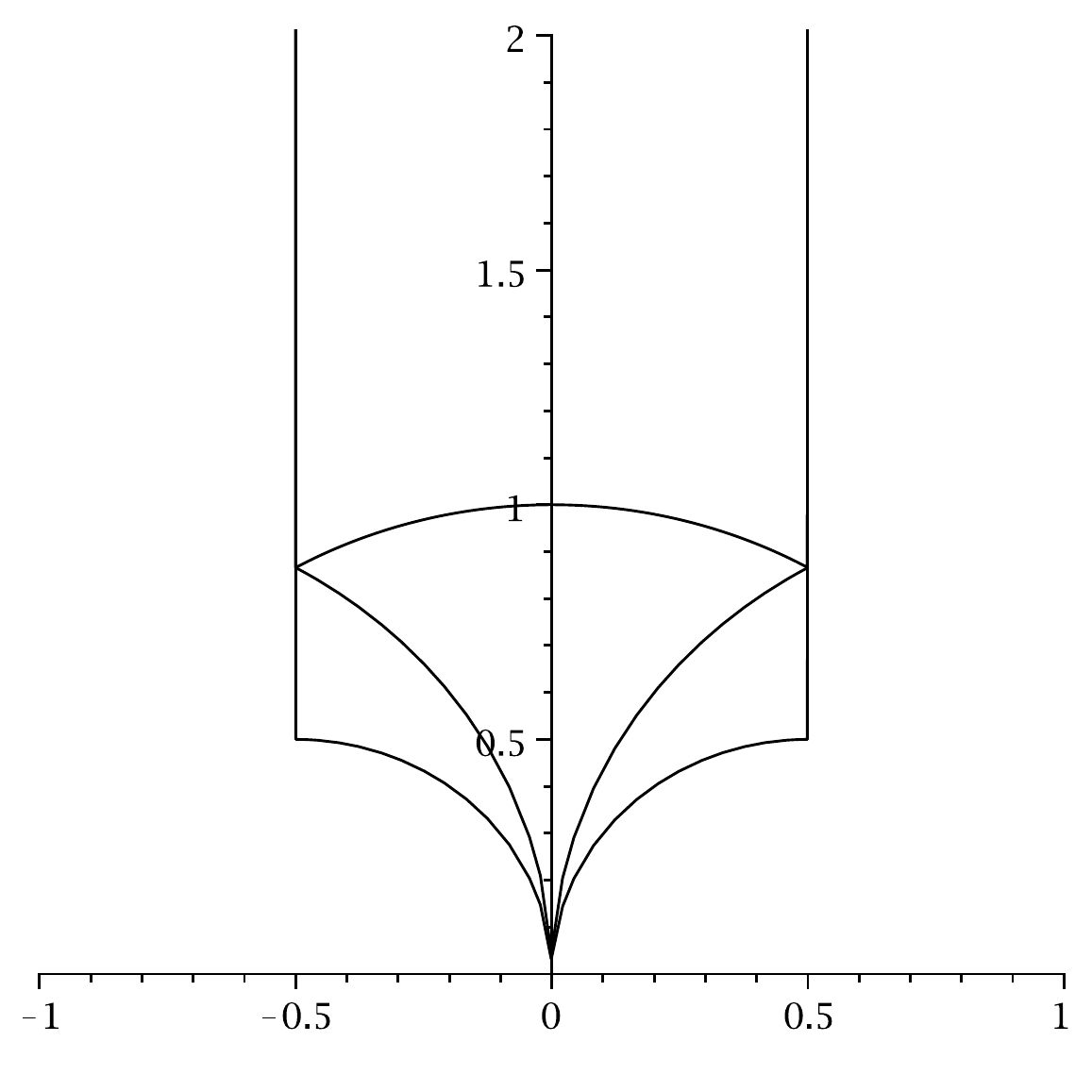}
\end{center}

Modular forms in $M_k^!(2)$ with real coefficients demonstrate a
nice property on the lower boundary of this fundamental domain.  For
a Fourier series $f(z) = \sum a(n) e^{2\pi i n z}$ with real Fourier
coefficients $a(n)$, note that $\overline{f(a+bi)} = \sum
\overline{a(n)} e^{-2\pi i n a} e^{-2\pi n b} = f(-a+bi)$, or
$\overline{f(z)} = f(-\overline{z})$. Thus, for modular forms $f$ of
weight $k$ on $\Gamma_0(2)$, if we let $z = -\frac{1}{2} +
\frac{1}{2}e^{i\theta}$ for $0 \leq \theta \leq \frac{\pi}{2}$, so
that $z$ is on the lower boundary of the (symmetric) fundamental
domain, we find that $\smat{1}{0}{2}{1}z = \frac{z}{2z+1} =
\frac{1}{2} - \frac{1}{2}e^{-i\theta} = -\overline{z}$.  Therefore,
$f(z) = (2z+1)^{-k} f\left(\frac{z}{2z+1}\right) = e^{-ik\theta}
\overline{f(z)}$, and the normalized modular form
$e^{\frac{ik\theta}{2}} f(-\frac{1}{2} + \frac{1}{2}e^{i\theta})$ is
real valued for $\theta$ between $0$ and $\frac{\pi}{2}$.

It is useful to define three particular modular forms of level 2. As
usual, let $q = e^{2\pi i z}$. Let
\[\psi(z) = \left(\frac{\eta(z)}{\eta(2z)}\right)^{24} = q^{-1} - 24
+ 276 q + \ldots \in M_0^\sharp(2) \] be the Hauptmodul for
$\Gamma_0(2)$. This form has integer coefficients, has a pole at
$\infty$, and vanishes at $0$. Moreover, by the above argument
$\psi(z)$ is real-valued on the lower boundary of the fundamental
domain, taking on values in $[-64,0]$. The special value
$\psi\left(-\frac{1}{2}+\frac{i}{2}\right) = -64$ arises from the
relationships between $\psi(z)$ and $j(z)$ given by
\[
j(z) = \frac{(\psi(z)+256)^3}{\psi(z)^2}, \qquad j(2z) =
\frac{(\psi(z)+16)^3}{\psi(z)}.
\]
By using the fact that $j\left(-\frac{1}{2}+\frac{i}{2}\right) =
j(i) = 1728$, it is easy to see that
$\psi\left(-\frac{1}{2}+\frac{i}{2}\right)$ must be a common root of
the polynomials $(x+256)^3-1728x^2 = (x+64)(x-512)^2$ and
$(x+16)^3-1728x = (x-8)^2(x+64)$.

Next, let
\[F_2(z) = 2E_2(2z) - E_2(z) = 1 + 24 \sum_n \left(\sum_{d|n, d \,
\rm{odd}} d \right) q^n\] be the unique normalized holomorphic
modular form of weight $2$ and level $2$.  Here $E_2(z)$ is the
weight 2 Eisenstein series $E_2(z) = 1-24 \sum_{n=1}^\infty
\sigma(n) q^n$. The form $F_2(z)$ has integer coefficients and a
single zero at the elliptic point $-\frac{1}{2} + \frac{i}{2}$. This
can be seen by noting that $F_2(\frac{z}{2z+1}) = F_2(z+1) = F_2(z)$
at $z = -\frac{1}{2} + \frac{i}{2}$, so an application of the
modular equation yields $-F_2(z) = F_2(z)$. Uniqueness comes from
applying the valence formula for $\Gamma_0(2)$ found in, for
instance,~\cite{El-G2}. Additionally, we define the Eisenstein
series $S_4(z) \in M_4(2)$ as
\begin{equation}
\label{S4}S_4(z) = \frac{E_4(z) - E_4(2z)}{240} =
q+8q^2+28q^3+\cdots;
\end{equation} it is easily checked that $S_4$ has
integral Fourier coefficients and vanishes at $\infty$.  It does not
vanish at the cusp at 0, as the valence formula shows that there are
no cusp forms of weight $4$ and level $2$.

We now use these forms to construct a basis for $M_k^\sharp(2)$. For
general even weight $k$, we write $k = 4\ell + k'$, where $k' \in
\{0, 2\}$.  A basis for $M_k^\sharp(2)$ is given by
\[f_{k, n}^{(2)}(z) = q^{-n} + O(q^{\ell+1}),\] for all integers
$n \geq -\ell$.  We note that this is an extension of the basis
given in~\cite{DoudJ} for $M_k(2)$; similar sequences of modular
forms for many levels appear in~\cite{El-G}. The basis elements are
constructed as follows.  We first define $F_0 = 1$ and set $f_{k,
-\ell}^{(2)}(z) = S_4^\ell(z) F_{k'}(z)$.  Next, for each $n >
-\ell$ we define $f_{k, n}^{(2)}(z)$ inductively by multiplying
$f_{k, n-1}^{(2)}(z)$ by $\psi(z)$ and subtracting off earlier basis
elements.   Note that since $S_4$, $F_{k'}$, and $\psi$ have
integral Fourier coefficients, each $f_{k, n}^{(2)}(z)$ has integral
Fourier coefficients and is of the form $S_4^\ell(z) F_{k'}(z)
F(\psi(z))$, where $F(x)$ is a polynomial with integer coefficients
of degree $n + \ell = n + \lfloor \frac{k}{4} \rfloor$. (In fact,
$F(x)$ is a generalized Faber polynomial; see~\cite{Fa1},
\cite{Fa2}.) Thus, if all of the zeros of the basis elements $f_{k,
n}^{(2)}$ lie on the lower boundary of the fundamental domain, then
all of the zeros of the polynomial $F(x)$ must lie in the interval
$[-64, 0]$.

The main result of this paper is the following theorem.
\begin{thm} \label{mainthm}
Let $f_{k, n}^{(2)}(z)$ be as above.  If $\ell \geq 0$ and $n \geq
14\ell+8$, or if $\ell <0$ and $n \geq 15|\ell|+8$, then at least
$\lfloor \frac{\sqrt{3}}{2}n + \frac{k}{6} \rfloor$ of the
$n+\lfloor \frac{k}{4} \rfloor$ nontrivial zeros of $f_{k,
n}^{(2)}(z)$ in the fundamental domain for $\Gamma_0(2)$ lie on the
lower boundary of the fundamental domain.
\end{thm}

We note that the bounds  $n \geq 14\ell+8$ and $n \geq 15|\ell|+8$
are not sharp, and that often many more of the zeros lie on the arc.
In fact, for certain weights $k$ close to $0$, all of the zeros of
all of the $f_{k, n}^{(2)}(z)$ are on the lower boundary of the
fundamental domain.  However, some restriction on $n$ in relation to
$\ell$ is necessary, as there are also examples of $f_{k,
n}^{(2)}(z)$ with zeros elsewhere.  We discuss this further in
Section~\ref{gamma03}.  Additionally, we obtain similar results for
a family of modular forms in $M_k^\sharp(2)$ whose coefficients are
dual to the Fourier coefficients of $f_{k, n}^{(2)}(z)$ and for a
basis for the space $M_k^\sharp(3)$, showing that many of the zeros
of these modular forms in the appropriate fundamental domain lie on
the lower boundary.

The remainder of this paper proceeds as follows: in
Section~\ref{genfn}, we give a generating function for the basis
elements $f_{k, n}^{(2)}(z)$ and approximate their values on the
lower boundary of the fundamental domain by a trigonometric
function.  In Section~\ref{bounds} we bound the error term to show
that most of the zeros lie on the appropriate arc.
Section~\ref{details} gives technical details on bounds for the
error, and in Section~\ref{gamma03} we discuss extensions of the
main theorem to other modular forms for $\Gamma_0(2)$ and to forms
for $\Gamma_0(3)$.

\section{Generating functions and integration} \label{genfn}
In this section we use Cauchy's integral formula to relate the basis
elements $f_{k, n}^{(2)}(z)$ to a trigonometric function.  Letting
$r = e^{2\pi i \tau}$, a generating function for the basis elements
$f_{k, n}^{(2)}(z)$ is given in El-Guindy's paper~\cite[Theorem
1.2]{El-G} as
\begin{equation}
\label{GenF} \sum_{n=-\ell}^\infty f_{k, n}^{(2)}(z) r^n =
\frac{(S_4^\ell F_{k'})(z)}{(S_4^\ell  F_{k'})(\tau)}
\frac{\psi(\tau) F_2(\tau)}{\psi(\tau)-\psi(z)}.\end{equation}

Multiply by $r^{-n-1}$ and integrate around $r=0$.  Changing
variables from $r$ to $\tau$ and noting that $\psi(\tau)F_2(\tau) = \frac{d}{d\tau} \psi(\tau)$, we obtain
\[f_{k, n}^{(2)}(z) = \int_{-\frac{1}{2} + i A}^{\frac{1}{2} + i A} \frac{(S_4^\ell
F_{k'})(z)}{(S_4^\ell F_{k'})(\tau)} \frac{\frac{d}{d\tau}
(\psi(\tau)-\psi(z))}{\psi(\tau)-\psi(z)} \frac{e^{-2\pi i n
\tau}}{-2\pi i} d\tau,\] where $A$ is some real number larger than
$\frac{1}{2}$. We move the contour downward, noting that there is a
pole whenever $\psi(\tau) - \psi(z)$ is zero, which happens when
$\tau$ is equivalent to $z$ under the action of $\Gamma_0(2)$. There
are no other poles, since $S_4$ has no zeros in the upper half
plane, and if $k' = 2$, the zero of $F_2(\tau)$ is canceled by the
zero of $\frac{d}{d\tau} \psi(\tau) = \psi(\tau) F_2(\tau)$. The
closed contour that is the difference between the old integral and
the new integral is in the clockwise direction, so we get a factor
of $-1$ in Cauchy's integral theorem, and at each pole, we obtain a
term of $-2\pi i$ times the residue of the integrand.

If a function $f(\tau)$ has a zero of order $n$ at $\tau_0$, its
logarithmic derivative has a simple pole with residue $n$.  In
calculating the residue, note that part of the integrand is the
logarithmic derivative of $\psi(\tau) - \psi(z)$, which has a simple
zero exactly at the values we are looking at, since $\psi$ is a
Hauptmodul for $\Gamma_0(2)$.  This means that the logarithmic
derivative of $\psi(\tau) - \psi(z)$ at points equivalent to $z$
under $\Gamma_0(2)$ will just give us a factor of $1$ in the
residue. Supposing that $\tau = \gamma z = \frac{az+b}{cz+d}$ for
some $\gamma \in \Gamma_0(2)$, this is then multiplied by
\[\frac{e^{-2\pi i n \gamma z}}{-2\pi i} \frac{(S_4^\ell F_{k'})(z)}
{(S_4^\ell F_{k'})(\gamma z)}.\] Since the denominator is a modular
form of weight $k$ on $\Gamma_0(2)$, we have $ (S_4^\ell
F_{k'})(\gamma z) = (cz+d)^k (S_4^\ell F_{k'})(z)$ for $\gamma =
\smat{a}{b}{c}{d}$, and the residue becomes $\frac{1}{-2\pi
i}(cz+d)^{-k} e^{-2\pi i n(\gamma z)}$.

Assume that $z$ is on the lower boundary of the fundamental domain
for $\Gamma_0(2)$. The first two points where $\tau$ is equivalent
to $z$ through which the contour moves are $\tau = z$ and $\tau =
\frac{z}{2z+1}$. Calculating the residues, we find that
\[f_{k, n}^{(2)}(z) - e^{-2\pi i n z} - (2z+1)^{-k} e^{-2\pi i n
(\frac{z}{2z+1})} = \int_C \frac{(S_4^\ell F_{k'})(z)}{(S_4^\ell
F_{k'})(\tau)} \frac{\frac{d}{d\tau} (\psi(\tau)-\psi(z))}
{\psi(\tau)-\psi(z)} \frac{e^{-2\pi i n \tau}}{-2\pi i} d\tau,\]
where $C$ is a contour that moves from left to right across the
fundamental domain and passes below the points $\tau = z$ and $\tau
=\frac{z}{2z+1}$ and above all other points equivalent to $z$ under
the action of $\Gamma_0(2)$.

We write $z = -\frac{1}{2} + \frac{1}{2} e^{i\theta}$ for some $0
\leq \theta \leq \frac{\pi}{2}$, so that $\frac{z}{2z+1} =
\frac{1}{2} -\frac{1}{2}e^{-i\theta}$; then the quantity $e^{-2\pi i
n z} + (2z+1)^{-k} e^{-2\pi i n \left(\frac{z}{2z+1}\right)}$ can be
simplified to \[(-1)^n e^{-i k \theta/2} e^{\pi n \sin \theta} 2
\cos\left(\frac{k\theta}{2} - \pi n \cos \theta\right).\]

Putting all of this together and multiplying through by
$e^{ik\theta/2} e^{-\pi n \sin \theta}$, we end up with
\[e^{ik\theta/2} e^{-\pi n \sin \theta} f_{k,
n}^{(2)}\left(-\frac{1}{2} + \frac{1}{2} e^{i\theta}\right) - (-1)^n
2 \cos\left(\frac{k\theta}{2} - \pi n \cos \theta\right) =
\]\[e^{ik\theta/2} e^{-\pi n \sin \theta} \int_C
\frac{(S_4^\ell F_{k'})(z)}{(S_4^\ell F_{k'})(\tau)}
\frac{\psi(\tau) F_2(\tau)}{\psi(\tau)-\psi(z)} e^{-2\pi i n \tau}
d\tau.\]  By the argument in Section~\ref{defs}, the left hand side
is a real-valued function of $\theta$. We note that the cosine
function takes on alternating values of $\pm 2$ whenever
$\frac{k\theta}{2} - \pi n \cos \theta$ is equal to $m\pi$ for $m\in
\Z$. Since this quantity moves from $-n\pi$ at $\theta = 0$ to
$\frac{k\pi}{4}$ at $\theta = \frac{\pi}{2}$, we know that there must be at
least $n + 1 + \lfloor \frac{k}{4} \rfloor$ times that this happens.
Thus, if we bound the integral term in absolute value by $2$, then
by the Intermediate Value Theorem we must have at least $n + \lfloor
\frac{k}{4} \rfloor$ zeros of the modular form on this arc.

The dimension of the space of holomorphic modular forms $M_k(2)$ is
$\lfloor \frac{k}{4} \rfloor + 1$, and we get at most $\lfloor
\frac{k}{4} \rfloor$ zeros not at elliptic points for these forms by
the general valence formula in~\cite{El-G2}. (Note that $F_2(z)$ has
a zero at $-\frac{1}{2} + \frac{i}{2}$, which is already on the arc
in question, so $f_{k, n}^{(2)}$ has a trivial zero there if $k' =
2$.) Counting a pole of order $n$ at $\infty$ and no other poles
gives us a total of $n + \lfloor \frac{k}{4} \rfloor$ zeros of the
basis element $f_{k, n}^{(2)}(z)\in M_k^\sharp(2)$ whose locations
are unknown.  This argument proves that if the weighted modular form
is close enough to the cosine function, then all of these zeros must
be simple and must be on this arc on the lower boundary of the
fundamental domain.

Unfortunately, it is difficult to move the contour down far enough
to prove that all of the zeros are on this arc; as the weight $k$ or
the order $n$ of the pole increases, the contour will need to get
closer and closer to $\tau = 0$ if $\theta$ is close to $0$. Recall
that for every fixed value of $z$, we need to prove that the
integral is bounded by 2, after moving the contour below that value
of $z$. As $z$ gets close to the real line, this becomes very
difficult--either the contour is not straight and the integral is
harder to estimate, as $\tau$ does not have a fixed imaginary part,
or the contour must pass through more residues, adding additional
terms to the equation.

We can still prove that the majority of the zeros do indeed lie on
this arc by choosing a fixed height for the contour, estimating the
value of the integral along that contour, and showing that its
absolute value is bounded above by $2$.  The goal then is to choose
a contour low enough to capture as many zeros as possible, yet high
enough to avoid additional residues and to avoid large values inside
the integral.  We choose $\tau = u + \frac{i}{5}$ for $\abs{u} \leq
\frac{1}{2}$, so that the contour has constant imaginary part
$\frac{1}{5}$.

For this choice to work, we must also limit the range of $z$, so
that our contour passes below $\tau = z$ and $\tau = \frac{z}{2z+1}$
but above other images of $z$ under $\Gamma_0(2)$.  If $z =
-\frac{1}{2} + \frac{1}{2} e^{i\theta}$ for $\frac{\pi}{6} \leq
\theta \leq \frac{\pi}{2}$, then $z$ has imaginary part $\geq
\frac{1}{4}$, and a contour at a height of $\frac{1}{5}$ picks up
residues at $\tau = z$ and $\tau = \frac{z}{2z+1}$ but no other
points equivalent to $z$ under the action of $\Gamma_0(2)$; the
maximum possible imaginary part of such a point is $\frac{1}{6}$. In
this case, the quantity $\left(\frac{k\theta}{2} - \pi n \cos
\theta\right)$ inside the cosine function has the value
$\frac{k\pi}{12} - \pi n \frac{\sqrt{3}}{2}$ at $\theta =
\frac{\pi}{6}$ and the value $\frac{k\pi}{4}$ at $\theta =
\frac{\pi}{2}$, and passes through at least $\lfloor \frac{k}{6} +
n\frac{\sqrt{3}}{2} \rfloor$ multiples of $\pi$. Bounding the
integral by $2$ for the appropriate $f_{k, n}^{(2)}(z)$ will finish
the proof of Theorem~\ref{mainthm}.

\section{Bounding the integral} \label{bounds}
In this section we bound \[e^{ik\theta/2} e^{-\pi n \sin \theta}
\int_C \frac{(S_4^\ell F_{k'})(z)}{(S_4^\ell F_{k'})(\tau)}
\frac{\psi(\tau) F_2(\tau)}{\psi(\tau)-\psi(z)} e^{-2\pi i n \tau}
d\tau\] for the values $z=-\frac{1}{2} + \frac{1}{2}e^{i\theta}$
with $\theta \in [\frac{\pi}{6}, \frac{\pi}{2}]$ and $\tau =
u+\frac{i}{5}$ with $u \in [-\frac{1}{2}, \frac{1}{2}]$. We will
also give some indication of how this bound might change if we allow
$\theta$ to approach $0$ and alter the countour accordingly. Details
for the computation of the numerical bounds that appear here are
provided in the next section.

We seek a bound for \[e^{\frac{ik\theta}{2}} e^{-\pi n \sin\theta}
f_{k, n}^{(2)}\left(-\frac{1}{2} + \frac{1}{2} e^{i\theta}\right) -
(-1)^n 2 \cos\left(\frac{k\theta}{2} - \pi n \cos \theta\right) =
\]\[e^{\frac{ik\theta}{2}} e^{-\pi n \sin\theta} \int_{-\frac{1}{2}}^
{\frac{1}{2}} \frac{(S_4^\ell F_{k'})(-\frac{1}{2} + \frac{1}{2}
e^{i\theta})}{(S_4^\ell F_{k'})(u+\frac{i}{5})}
\frac{\psi(u+\frac{i}{5}) F_2(u+\frac{i}{5})}
{\psi(u+\frac{i}{5})-\psi(-\frac{1}{2} + \frac{1}{2} e^{i\theta})}
e^{-2\pi i n u} e^{\frac{2\pi n}{5}} du,\] which is real-valued, by
something less than $2$. In absolute value, this integral is
\[e^{-\pi n (\sin \theta - \frac{2}{5})} \left|
\int_{-\frac{1}{2}}^{\frac{1}{2}} \frac{(S_4^\ell
F_{k'})(-\frac{1}{2} + \frac{1}{2} e^{i\theta})}{(S_4^\ell
F_{k'})(u+\frac{i}{5})} \frac{\psi(u+\frac{i}{5})
F_2(u+\frac{i}{5})} {\psi(u+\frac{i}{5})-\psi(-\frac{1}{2} +
\frac{1}{2} e^{i\theta})} e^{-2\pi i n u} du\right|.\]

Consider the exponential term $e^{-\pi n (\sin \theta -
\frac{2}{5})}$.  We have chosen $\Im(\tau) = \frac{1}{5}$, so that
$\sin\theta - \frac{2}{5} >0$ for $\theta \in [\frac{\pi}{6},
\frac{\pi}{2}]$, and this term has exponential decay as $n
\rightarrow \infty$; in this case $e^{-\pi (\sin \theta -
\frac{2}{5})} < .73041$.  If we find an upper bound for the absolute
value of the integral, then for large enough $n$ the right hand side
is indeed less than 2, and we can apply the Intermediate Value
Theorem as desired.  It turns out that we can find a bound for the
absolute value that removes the dependence on $n$, but may be
exponential in $\ell$. However, if $n$ is large enough in relation
to $\ell$, then we will see that for a fixed weight $k$, all but
finitely many of the $f_{k, n}^{(2)}(z)$ have zeros on the
appropriate arc.

We note that the absolute value of the integral is certainly bounded
above by \[ \int_{-\frac{1}{2}}^{\frac{1}{2}}
\left|\frac{S_4(-\frac{1}{2} + \frac{1}{2}
e^{i\theta})}{S_4(u+\frac{i}{5})}\right|^\ell \left|
\frac{F_{k'}(-\frac{1}{2} + \frac{1}{2} e^{i\theta})
F_2(u+\frac{i}{5}))} {F_{k'}(u+\frac{i}{5})}\right|
\left|\frac{\psi(u+\frac{i}{5})} {\psi(u+\frac{i}{5})
-\psi(-\frac{1}{2} + \frac{1}{2} e^{i\theta})}\right| du,\] and
already the dependence on $n$ has vanished.  For the terms involving
$S_4$ and $F_{k'}$, we find an upper bound for the maximum possible
value of these terms over the appropriate ranges of $u$ and
$\theta$, and pull these upper bounds outside of the integral. This
leaves us with the contribution from
\[\int_{-\frac{1}{2}}^{\frac{1}{2}} \left|\frac{\psi(u+\frac{i}{5})}
{\psi(u+\frac{i}{5}) -\psi(-\frac{1}{2} + \frac{1}{2}
e^{i\theta})}\right| du,\] which we consider in pieces.

Consider first the quantity \[\left| \frac{S_4(-\frac{1}{2} +
\frac{1}{2} e^{i\theta})}{S_4(u+\frac{i}{5})}\right|^\ell.\]
Computations, explained in more detail in the next section, yield
\[.014 \leq \left|S_4\left(u+\frac{i}{5}\right)\right|\leq
2.44141.\] As a power series in $q$, we know that $S_4$ has positive
coefficients, as seen in \eqref{S4}. The maximum value of
$|S_4(\tau)|$ occurs when $u=0$, when $q$ is real and positive.
Heuristically, the minimum value should occur when there is maximum
cancelation between terms, or when $q$ is real and negative, so that
$u = \pm \frac{1}{2}$, and we confirm this computationally.
Decreasing $\Im(\tau)$ both increases the upper bound, as we are
adding larger positive terms, and potentially decreases the lower
bound due to cancelation.

Heuristically, the maximum value of $|S_4(z)|$ should occur when
$(-1+\cos\theta)$ is close to $0$, meaning $\theta$ is close to $0$,
as here $q$ is real, positive, and close to 1.  Similarly, the
minimum value of $S_4(z)$ should occur when $(-1+\cos\theta)$ is
close to $\pm 1$, meaning $\theta$ is close to $\frac{\pi}{2}$; here
$q$ is real but negative, so there is extensive cancelation when
adding terms.  In this case, though, the size of $q$ depends on
$\theta$, as we have $|e^{2\pi i \left(-\frac{1}{2} +
\frac{1}{2}e^{i\theta}\right)}|= |e^{-\pi\sin(\theta)}|$.
Computationally, the minimum indeed occurs at $\theta =
\frac{\pi}{2}$, where $|e^{-\pi\sin(\frac{\pi}{2})}| \approx
.04322$.  In general, we have \[|e^{2\pi i \left(-\frac{1}{2}
+\frac{1}{2}e^{i\theta}\right)}|= |e^{-\pi\sin(\theta)}|\leq
e^{-\pi/2} \approx .20788.\]   We compute that \[.03 \leq
\left|S_4\left(-\frac{1}{2} + \frac{1}{2} e^{i\theta}\right)\right|
\leq .99995.\] Moving the lower bound on $\theta$ closer to $0$
increases the maximum value of $|S_4(z)|$, though it does not appear
to affect the minimum value.

Putting this together, we have, for $\ell \geq 0$,
\[\left|\frac{S_4(-\frac{1}{2} + \frac{1}{2} e^{i\theta})}
{S_4(u +\frac{i}{5})}\right|^\ell \leq |71.425|^\ell,\] and for
$\ell < 0$, \[\left|\frac{S_4(-\frac{1}{2} + \frac{1}{2}
e^{i\theta})} {S_4(u +\frac{i}{5})}\right|^\ell \leq
|81.38034|^{\abs{\ell}}.\]

Next, we consider the term
\[\left|\frac{F_{k'}(-\frac{1}{2} + \frac{1}{2} e^{i\theta})
F_2(u+\frac{i}{5})}{F_{k'}(u + \frac{i}{5})}\right|. \] If $k' = 2$,
this is $\left|F_2(-\frac{1}{2} + \frac{1}{2} e^{i\theta})\right|$,
which is bounded above by 8.00067.  If $k' = 0$, this is
$\left|F_2(u+\frac{i}{5})\right|$, which is bounded above by
12.50005. Either way, the contribution is no more than 12.50005.
Note that $F_2(\tau)$ has positive coefficients, and is therefore
large when $S_4(\tau)$ is large.

Finally, we consider
\[ \int_{-\frac{1}{2}}^{\frac{1}{2}}\left|\frac{\psi(u+\frac{i}{5})}
{\psi(u+\frac{i}{5})-\psi(-\frac{1}{2} + \frac{1}{2}
e^{i\theta})}\right|~du.\]  The Hauptmodul $\psi(z)$ is an injective
mapping on the fundamental domain; it is real-valued and strictly
decreasing on $\theta \in [0,\frac{\pi}{2}]$. We restrict ourselves
to $\theta \in [\frac{\pi}{6},\frac{\pi}{2}]$. Computation shows
that $\psi(z)\in [-64, -.033]$ on this domain.  In contrast,
$\psi(u+\frac{i}{5})$ takes on a wide range of values, including
some with very large and some with very small modulus. Bounding the
numerator and denominator separately yields a trivial upper bound of
roughly 79000, while numerical calculations indicate that the actual
maximum is a little larger than 1.

In order to achieve a sharper bound, we will instead consider a
related quantity,
\begin{equation}
\label{D}
D(z,\tau) =
\frac{\psi(z)}{\psi(\tau)-\psi(z)},
\end{equation}
where we use $\psi(z)$ to indicate $\psi(-\frac{1}{2} + \frac{1}{2}
e^{i\theta})$ and $\psi(\tau)$ to indicate $\psi(u+\frac{i}{5})$ for
ease of notation.  This quantity is related to our Hauptmodul
expression by the identity
\begin{equation}
\label{relation} \left|\frac{\psi(\tau)}{\psi(\tau)-\psi(z)}\right|
= \left|1+\frac{\psi(z)} {\psi(\tau)-\psi(z)}\right| =
\left|1+D(z,\tau)\right|.
\end{equation}
It is easier to work with $D(z,\tau)$, as we know the numerator is
real valued and within a small range, and bounding this quantity
proves useful in Section~\ref{gamma03} when discussing extensions of
Theorem~\ref{mainthm}.

We will break our path of integration into pieces, and consider
$\psi(\tau)$ in relation to $\psi(z)$ on each.  As we know that
$\psi(z)$ is real, we consider the real and imaginary parts of
$\psi(\tau) = \psi\left(u+\frac{i}{5}\right)$ separately for $u \in
[-\frac{1}{2}, \frac{1}{2}]$. It is clear that $e^{2\pi i \tau} =
e^{-\frac{2\pi}{5}}\left(\cos(2\pi u)+i\sin(2\pi u)\right)$, and so
$\Re(\psi(u+\frac{i}{5})) = \Re(\psi(-u+\frac{i}{5}))$, while
$\Im(\psi(u+\frac{i}{5})) = -\Im(\psi(-u+\frac{i}{5}))$.  With this
in mind, we restrict our calculations to $u \in [-\frac{1}{2}, 0]$
and use symmetry for $u \in [0, \frac{1}{2}]$.

The numerical techniques described in the next section reveal that
on the interval $u \in [-.5, -.21516]$, we have either
$\Re(\psi(\tau)) > 0$, $\Re(\psi(\tau)) < -128$, or
$|\Im(\psi(\tau))| > 64$; it follows that $|D(z,\tau)| < 1$ on this
interval.

Next, we note that since $\psi(z)$ is real, then if we have the
bound $\Im(\psi(\tau)) \geq A > 0$, it follows that
\[
\left|D(z,\tau)\right|\leq
\left|\frac{\psi(z)}{(\Re(\psi(\tau))-\psi(z))+Ai}\right|.
\]
If $\Re(\psi(\tau)) \geq 0$, then this is bounded above by $1$,
while if $\Re(\psi(\tau)) < 0$, then the maximum possible value of
the right hand side for a fixed $\tau$ as $z$ varies is
\[
\sqrt{\left(\frac{\Re(\psi(\tau))}{A}\right)^2+1},
\]
occurring when $\psi(z) = \frac{A^2+\Re(\psi(\tau))^2}
{\Re(\psi(\tau))}$.  A lower bound for $\Re(\psi(\tau))$ thus gives
us an upper bound for $\abs{D(z, \tau)}$.

For $u \in [-.21516, -.18884]$ we have $\Im(\psi(\tau))
> 1$ and $\Re(\psi(\tau)) > -.0175.$  If $\Re(\psi(\tau))\geq 0$,
then our bound is 1, and if $0> \Re(\psi(\tau)) > -.0175$ then  our
bound is 1.00016. Either way, for $u \in [-.21516, -.18884]$, we
have $|D(z,\tau)|<1.00016$. Similarly, on $[-.18884, -.12878]$ we
have $\Im(\psi(\tau))\geq .033$, and we obtain a bound of 1.13192.

Finally, we consider $[-.12878, -0]$.  We have
$\Re(\psi(\tau))>-.01424$, and so
\[
\left|D(z,\tau)\right| \leq
\left|\frac{\psi(z)}{-.01424-\psi(z)}\right|.
\]
The maximum value occurs when $\psi(z) = -.033$, and we find $|D(z,
\tau)| < 1.75344$ here.

Altogether, we have $|D(z,\tau)| < 1.75344$.  By breaking the
integral into pieces, we compute more precisely that
\[
\int_{-\frac{1}{2}}^{\frac{1}{2}}\left|D(z,\tau)\right|~du \leq 2 \cdot 0.60496 =
1.20992.
\]
The relationship in \eqref{relation} allows us to conclude that
\[
\int_{-\frac{1}{2}}^{\frac{1}{2}}\left|\frac{\psi(u+\frac{i}{5})}
{\psi(u+\frac{i}{5})-\psi(-\frac{1}{2} + \frac{1}{2}
e^{i\theta})}\right| du  < 2.20992.
\]
We can improve this further by noting that the fact that
$\Re(\psi(\tau))
> 0$ on $[-.45787, -.22531]$ implies that the integrand is bounded
by 1 here, yielding
\[
\int_{-\frac{1}{2}}^{\frac{1}{2}}\left|\frac{\psi(u+\frac{i}{5})}
{\psi(u+\frac{i}{5})-\psi(-\frac{1}{2} + \frac{1}{2} e^{i\theta})}
du \right| < 1.74520.
\]
We see that if we extend $\theta$ closer to $0$, and hence also
decrease $\Im(\tau)$, this term has the most potential to blow up
near $u=0$, as this is where $\psi(z)$ and $\Re(\psi(\tau))$ are
both small.  Additionally, if $0 \leq \theta \leq \frac{\pi}{2}$,
then any contour with fixed imaginary part less than $\frac{1}{2}$,
such as $\tau = u+\frac{i}{5}$, will cross the arc $z=-\frac{1}{2} +
\frac{1}{2}e^{i\theta}$ and at least one of its images under
$\Gamma_0(2)$, so restricting our $\theta$ values is necessary to
avoid a zero in the denominator.

Putting all of these pieces together and using the fact that
$\sin\theta$ is decreasing on $[\frac{\pi}{6}, \frac{\pi}{2}]$, we
see that for $\ell
\geq 0$,
\begin{align*}&e^{-\pi n (\sin \theta - \frac{2}{5})} \left|
\int_{-\frac{1}{2}}^{\frac{1}{2}} \frac{(S_4^\ell
F_{k'})(-\frac{1}{2} + \frac{1}{2} e^{i\theta})}{(S_4^\ell
F_{k'})(u+\frac{i}{5})} \frac{\psi(u+\frac{i}{5})
F_2(u+\frac{i}{5})} {\psi(u+\frac{i}{5})-\psi(-\frac{1}{2}
+ \frac{1}{2} e^{i\theta})} e^{-2\pi i n u} du\right|\\
& < .73041^n|71.425|^\ell (12.50005)(1.74520).
\end{align*}
Note that $(.73041^n) (12.50005)(1.74520) < 2$ if $n\geq 8$, and
$(.73041^n)71.425 < 1$ if $n\geq 14$; hence, the integral is less
than our desired bound $2$ if $\ell \geq 0$ and $n \geq 14\ell+8$.
Similarly, for $\ell < 0$, we replace $|71.425|^\ell$ with
$|81.38034|^{\abs{\ell}}$, and find that our integral is bounded by
$2$ if $n \geq 15\abs{\ell}+8$.  We can then apply the Intermediate
Value Theorem to prove that the appropriate number of zeros are on
the desired arc.

\section{Rigorously computing upper and lower bounds}
\label{details}

In the previous section, while bounding our integral we used upper
and lower bounds on Eisenstein series and the Hauptmodul for values
on a circular arc on the boundary of the fundamental domain and on a
straight line segment.  In this section we justify those bounds.

It is useful for most of these calculations to truncate each series.
For a modular form $f$ with Fourier series $f = \sum a_f(n) q^n$, we
will choose a positive integer $N$ and let $\tilde{f}$ be the
truncation of the Fourier series of $f$ up to and including the
$q^{N}$ term, and we let $Rf = f - \tilde{f}$ be the remaining tail
of the series. We bound $\tilde{f}$ and $Rf$ separately.

The calculations for the Eisenstein series are straightforward, as
we have explicit formulas for the Fourier coefficients in terms of
divisor functions, while calculations for $\psi$ require a little
more finesse. We do not have a nice formula or a sharp bound for the
growth rate of the Hauptmodul coefficients, and they are quite
large, so there are more terms making a significant contribution to
the value of the series. We begin by bounding the values of the
Eisenstein series, and then use those bounds to tame $\psi(z)$.

By \eqref{S4} we have
\[
S_4(z)=\sum_{n=1}^\infty \left(\sigma_3(n)-\sigma_3 \left(
\frac{n}{2} \right) \right)q^n.
\]
For $k\geq 1$, we can generously bound $\sigma_k(n) = \sum_{d|n}d^k$
by $\sqrt{n}\cdot\sqrt{n}^{k}+\sqrt{n}\cdot n^k$ by considering
pairs of divisors $\left( d, \frac{n}{d} \right)$.  Thus,
$\sigma_k(n)\leq n^{\frac{k+1}{2}}+n^{k+1}$.  If $|e^{2\pi i z}|\leq
t$, then we can bound $RS_4(z)$ by
\begin{align*}
\left|\sum_{n=N+1}^\infty \left(\sigma_3(n)-\sigma_3\left(
\frac{n}{2}\right)\right) q^{n} \right| \leq & \sum_{n=N+1}^{n=N+10}
\left(\sigma_3(n)-\sigma_3\left( \frac{n}{2} \right) \right) t^{n} + \\
& \sum_{n=N+11}^\infty (n^2+n^4)t^{n} -
\sum_{n=\left\lfloor{\frac{N+11}{2}} \right\rfloor}^\infty
(1+n^3)t^{2n}.
\end{align*}
The last term reflects a lower bound on the summand $\sigma_3\left(
\frac{n}{2}\right)$ coming from the even terms.  Standard Taylor
series methods involving derivatives of the geometric series
$(1-x)^{-1} = \sum x^n$ taken at $x=t$ allow us to bound each
infinite series.  For example,
\[
\sum_{n=N+11}^\infty (n^2+n^4)t^{n} =  \frac{t^2+t}{(1-t)^3}+
\frac{t^4+11t^3+11t^2+t}{(1-t)^5} - \sum_{n=1}^{N+10}
(n^2+n^4)t^{n}.
\]
We take $N=50$. For $S_4(z)$, we note that $\abs{e^{2\pi i z}} \leq
t=e^{-2\pi\sin(\frac{\pi}{6})} =e^{\frac{-2\pi}{4}}$. We find that
$|RS_4(z)|\leq 2.86404 \cdot 10^{-23}$.

To find an upper bound for $|S_4(z)|$, we explicitly compute
$\sum_{n=1}^{50} \left(\sigma_3(n)-\sigma_3\left(
\frac{n}{2}\right)\right) t^n$ with $t=e^{\frac{-2\pi}{4}}$ and then
add the tail.  We find that $|S_4(z)|\leq .99995$ on the appropriate
arc.  For a lower bound on $|S_4(z)|$, we consider the real and
imaginary parts of $S_4(-\frac{1}{2} +\frac{1}{2}e^{i\theta})$
separately.  We compute the first derivative with respect to
$\theta$ for each part, and for $\theta\in [\frac{\pi}{6},
\frac{\pi}{2}]$, each derivative has a trivial upper bound of
$\sum_{n=1}^{50}\pi n \cdot \left( \sigma_3(n) -\sigma_3 \left(
\frac{n}{2} \right) \right) (e^{-2\pi/4})^n < 8.01$.  A computation
reveals that either the real or the imaginary part of $S_4(z)$ is at
least .0302 when $z = -\frac{1}{2}+\frac{1}{2}e^{i\theta}$ with
$\theta = \pi(\frac{1}{6}+\frac{1}{3}\cdot \frac{n}{40000})$, for
all $0 \leq n \leq 40000$.  The bounds on derivatives and on the
tail $RS_4(z)$ limit how close $\abs{S_4(z)}$ can get to $0$, and we
therefore conclude that $\abs{S_4(z)}$ is bounded below by $.03$.

We can do similar upper bound calculations for $S_4(\tau)$,
$F_2(\tau)$, and $F_2(z)$, using the additional fact that
$\abs{e^{2\pi i \tau}} \leq t=e^{-2\pi/5}$.  For $F_2$, we bound the
tail by
\[ \abs{RF_2(\tau) } =
\left|\sum_{n=N+1}^\infty a_f(n)q^{n}\right|\leq
\sum_{n=N+1}^{n=N+10} a_{F_2}(n)t^{n}+24\sum_{n=N+11}^\infty
\left(\frac{n}{2}+n^2\right) \left(e^{-\frac{2\pi}{5}}\right)^{n},
\]
where the use of $\frac{n}{2}$ instead of $n$ comes from the fact
that we only consider odd divisors.  We again take $N=50$ in each
case, and find that $|S_4(\tau)|\leq 2.44141$, $|F_2(z)| \leq
8.00067$, and $|F_2(\tau)| \leq 12.50005$.

To compute a lower bound for $|S_4(\tau)|$, we trivially bound the
derivatives of the real and imaginary parts of $S_4(u+\frac{i}{5})$
by \[\sum_{n=1}^{50}2\pi n \cdot a_{S_4}(n)(e^{-2\pi/5})^n <
48.83.\] We then compute the real and imaginary parts of
$S_4(u+\frac{i}{5})$ for $u = -\frac{1}{2}+\frac{1}{2}\cdot
\frac{n}{123000}$, where $0\leq n \leq 123000$, verifying that at
least one of these values is larger than .014010. We conclude that
$|S_4(\tau)|$ is bounded below by .014 on $u\in [-\frac{1}{2},0]$,
and use symmetry to extend this to $[0, \frac{1}{2}]$.

We now consider the Hauptmodul $\psi(\tau)$.  In the previous
section we needed information about the size of the real and
imaginary parts of $\psi(\tau)$ and the value of $\psi(z)$.  For
these computations, we work with the truncations $\tilde{\psi}(z)$
and $\tilde{\psi}(\tau)$, taking into account the growth of the real
and imaginary parts of the truncations and the error caused by
ignoring the tail.  In this case, the trivial bound on the partial
derivatives is much larger than in the Eisenstein series case, so we
truncate each series up to and including the $q^{30}$ term to
shorten our computation time.

We can express $\psi$ in terms of Eisenstein series of level 2 as
\[
\psi(z) = \frac{E_4(2z)}{S_4(z)}-16;
\]
we use this representation to bound $R\psi(z)$.  Observe that if we
truncate $\psi(z)$, then the tail satisfies
\begin{align*}
R\psi(z)&= \psi(z)-\tilde{\psi}(z)\\
&=\frac{E_4(2z)}{S_4(z)}-16-\tilde{\psi}(z)\\
&=\left(\tilde{\psi}(z)+16+\frac{E_4(2z)-(\tilde{\psi}(z)+16)
\tilde{S}_4(z)}{\tilde{S_4}(z)}\right)\frac{\tilde{S_4}(z)}{S_4(z)}
-16-\tilde{\psi}(z)\\
&\leq \left|\tilde{\psi}(z)+16\right|\left|\frac{\tilde{S}_4(z)}
{\tilde{S}_4(z)-RS_4(z)}-1\right|+ \frac{\left|\tilde{E_4}(2z)-
(\tilde{\psi}(z)+16)\tilde{S}_4(z)\right|+\left|RE_4(2z)\right|}
{\left|\tilde{S}_4(z)-RS_4(z)\right|}.
\end{align*}

We can now bound $R\psi(z)$.  We compute bounds for the $q^{31}$ to
$q^{40}$ terms directly to find that they contribute at most
$6.46551 \cdot 10^{-8}$, and use the truncation formula to bound the
remainder of the tail. In applying the formula, we truncate all
Eisenstein series at $N=50$ to use the previously computed upper and
lower bounds.  Note that, working as before, we can use the fact
that $E_4(2z) = 1 + 240\sum_{n=1}^{25} \sigma_3(n)q^{2n}$ to prove
that $|RE_4(2z)| \leq 6.40309 \cdot 10^{-29}$ and $|RE_4(2\tau)|
\leq 2.16794 \cdot 10^{-22}$.   We compute bounds for
$\tilde{\psi}(z)+16$ and $\tilde{E_4}(2z)- (\tilde{\psi}(z)+16)
\tilde{S}_4(z)$ by summing $\sum |a_f(n)|e^{\frac{-2\pi n}{4}}$ for
each series.  We find that the bounds are 544.01429 and $7.29909
\cdot 10^{-13}$, respectively. Putting all of this together we have,
for $N = 30$,
\begin{align*}
|R\psi(z)|\leq & \, 6.46551 \cdot 10^{-8}+544.01429 \cdot \left(
\frac{.03}{.03-7.05863 \cdot 10^{-30}}-1\right) + \\
& \frac{7.29909 \cdot 10^{-13}+ 6.40309 \cdot 10^{-29}}
{.03-7.05863 \cdot 10^{-30}} \\
< & \, 6.46754 \cdot 10^{-8}.
\end{align*}

We repeat the same calculations for $\psi(\tau)$.  We find that
\begin{align*}
|R\psi(\tau)| \leq & \, .001371+2593.07795\cdot \left(
\frac{.014}{.014-2.86404 \cdot 10^{-23}}-1\right) + \\
& \frac{6.40510\cdot 10^{-7}+ 2.16794 \cdot 10^{-22}}
{.014-2.86404 \cdot 10^{-23}} \\
< & \, 0.00142.
\end{align*}

Now that we have bounds on the error caused by truncation, we
numerically compute that $\psi(-\frac{1}{2} + \frac{1}{2}
e^{i\theta}) \leq -.033$ for $\theta \in [\frac{\pi}{6},
\frac{\pi}{2}]$ by calculating $\tilde{\psi}(z)$ for $N = 30$ at
$\theta = \frac{\pi}{6}$ and adding the bound for the tail
$|R\psi(z)|$.

The bounds for $\psi(\tau)$ are slightly more difficult. As with
lower bounds for the Eisenstein series, we consider the real and
imaginary parts of $\psi(\tau)$ separately. The maximum possible
growth rate for the real part is
\[
\left|\frac{d}{d u}\sum_{n=-1}^{30}a_\psi(n) e^{ \frac{-2\pi n}{5}}
\cos(2\pi n u) \right|\leq \sum_{n=-1}^{30} 2\pi n
\cdot|a_\psi(n)|\left(e^{\frac{-2\pi}{5}}\right)^n < 101197.78,
\]
and the same bound holds for the derivative of the imaginary part.
We again compute values of the real and imaginary parts of
$\psi(\tau)$ at a sampling of points, and use these bounds on the
derivatives to find intervals over which the real and imaginary
parts of $\psi(\tau)$ fall within the ranges given in the previous
section.

\section{Extensions of Theorem \texorpdfstring{\ref{mainthm}}{1}} \label{gamma03}

In this section we discuss the sharpness of Theorem \ref{mainthm},
extend Theorem \ref{mainthm} to a dual family of weakly holomorphic
modular forms for $\Gamma_0(2)$, and consider analogous theorems for
families of modular forms for other genus zero subgroups.

We have shown that if $f_{k,n}^{(2)}(z)$ is a basis element for the
space $M_k^\sharp(2)$ and $n$ is large enough compared to $|\ell|$,
then the majority of zeros of $f_{k, n}^{(2)}(z)$ lie on the lower
boundary of the fundamental domain for $\Gamma_0(2)$.  We note that
the bounds on $n$ are not sharp.  For example, if $\ell >0$, then
the $F_2(u+\frac{i}{5})$ term in the numerator of the integrand
takes on its largest values near $u=0$, while the
$S_4(u+\frac{i}{5})$ term in the denominator takes on its smallest
values near $u=\pm \frac{1}{2}$, and we have simply taken absolute
upper bounds for these functions on the interval $u \in
[-\frac{1}{2}, \frac{1}{2}]$ without accounting for interaction
between these terms. It is clear that sharper bounds on $n$ are
possible, and it is natural to ask if the zeros of $f_{k,
n}^{(2)}(z)$ always lie on $z=-\frac{1}{2}+\frac{1}{2}e^{i\theta}$.
In fact, we now exhibit explicit examples for which the zeros are
not on this arc.

Note that if $k' = 0$ and the degree $n + \ell = n + \frac{k}{4}$ of
the Faber polynomial $F(x)$ is equal to $1$, then we can directly
compute that $F(x) = x - (8\ell-24)$. For the root of this
polynomial to be in $[-64, 0]$, we must have $-5 \leq \ell \leq 3$.
Thus, if $\ell \geq 4$ or $\ell \leq -6$, the single nontrivial zero
of the modular form $f_{k, -\ell+1}(z)$ will not be on the lower
boundary of the fundamental domain.  If $k' = 2$, the analogous
polynomial is $x-8\ell$, and we must have $-8 \leq \ell \leq 0$ for
the zero to be on the lower boundary.  Similarly, if $k' = 0$ and
the degree $n + \ell$ of the polynomial $F(x)$ is equal to $2$, we
obtain the polynomial $x^2 + (48-8\ell) x + (32\ell^2-188\ell+24)$,
which has complex roots when $\ell < -6$ or $\ell \geq 6$.  Such
calculations can also be done for polynomials of other small degrees
to find ranges of $\ell$ and $n$ which will not work.  Thus, the
theorem is not true in general, and some condition on the size of
$n$ compared to $\ell$ is necessary.

Next, we consider a straightforward extension of
Theorem~\ref{mainthm} to obtain results for another family of
modular forms for $\Gamma_0(2)$.  Let $g_{k, n}^{(2)}(z)$, for $n
\geq -\ell+1$, be the unique modular form in $M_k^\sharp(2)$ which
vanishes at the cusp $0$ and has Fourier expansion beginning
\[g_{k, n}^{(2)}(z) = q^{-n} + O(q^\ell).\]
These $g_{k, n}^{(2)}(z)$ form a basis for the subspace of
$M_k^\sharp(2)$ consisting of forms that vanish at $0$.  It is
straightforward to construct $g_{k, n}^{(2)}(z)$ by setting $g_{k,
-\ell+1}^{(2)}(z) = S_4^\ell F_{k'} \psi(z)$ and obtaining the
$g_{k, n}^{(2)}$ for larger values of $n$ inductively by multiplying
earlier basis elements by powers of $\psi(z)$ and subtracting the
appropriate terms. These $g_{k, n}^{(2)}$ were studied in the $k=0$
case by Ahlgren~\cite{Ahlgren}, and generating functions for
arbitrary even weight $k$ were found by El-Guindy~\cite{El-G}.  We
note their similarity to the $f_{k, n}^{(2)}(z)$. In fact,
\[\sum_{n=-\ell+1}^\infty g_{k, n}^{(2)}(z) q^n =  \frac{(S_4^\ell
\psi F_{k'})(z)}{(S_4^\ell \psi F_{k'})(\tau)} \frac{\psi(\tau)
F_2(\tau)}{\psi(\tau)-\psi(z)} = -\sum_{m \geq \ell} f_{2-k,
m}^{(2)}(\tau) q^m,
\] and replacing $k$ with $2-k$ and switching $\tau$ and
$z$ shows that the $m$th Fourier coefficient of $f_{k, n}^{(2)}(z)$
is equal to the negative of the $n$th Fourier coefficient of
$g_{2-k, m}^{(2)}(z)$.

With the above generating function, we can repeat the calculations
in Section~\ref{genfn}, replacing $f_{k, n}^{(2)}$ with $g_{k,
n}^{(2)}$, to see that the quantity $e^{\frac{ik\theta}{2}} e^{-\pi
n \sin\theta} g_{k, n}^{(2)}(-\frac{1}{2}+ \frac{1}{2} e^{i\theta})$
is real valued, and that its difference from $(-1)^n 2
\cos\left(\frac{k\theta}{2} - \pi n \cos \theta\right)$ is given by
an integral almost identical to the integral for the $f_{k,
n}^{(2)}(z)$, except that the $\psi(\tau)$ in the numerator of the
integrand becomes a $\psi(z)$; hence our Hauptmodul expression is
the $D(z,\tau)$ function defined in~\eqref{D}. Following the
argument as before and using the bound for $\int \abs{D(z,
\tau)}~du$ in Section~\ref{bounds}, we find that if $\ell \geq 0$,
the difference is at most $(.73041)^n (71.425)^\ell (12.50005)
(1.20992)$, while if $\ell < 0$, it is at most $(.73041)^n
(81.38034)^{\abs{\ell}} (12.50005) (1.20992)$.  We obtain the
following theorem.
\begin{thm}
Let $g_{k, n}^{(2)}(z)$ be as above.  If $\ell \geq 0$ and $n \geq
14\ell+7$, or if $\ell <0$ and $n \geq 15|\ell|+7$, then at least
$\lfloor \frac{k}{6} + n\frac{\sqrt{3}}{2} \rfloor$ of the
$n+\lfloor \frac{k}{4} \rfloor - 1$ nontrivial zeros of $g_{k,
n}^{(2)}(z)$ in the fundamental domain for $\Gamma_0(2)$ lie on the
lower boundary of the fundamental domain.
\end{thm}

Finally, we consider extensions of this method to other subgroups.
Consider $\Gamma_0(N) = \{\smat{a}{b}{c}{d} \in \SL_2(\Z) : c \equiv
0 \pmod{N} \}$. If $N=2,3,5,7,13$, then this is a genus zero
subgroup with two inequivalent cusps and with Hauptmodul
\[
\psi_p(z)=\left(\frac{\eta(z)}{\eta(pz)}\right)^{\frac{24}{p-1}}.
\]
It is again possible to define analogous families $f_{k,
n}^{(N)}(z)$ and $g_{k, n}^{(N)}(z)$ for these groups, and to
use~\cite{El-G} to find a generating function for each family that
mirrors \eqref{GenF}.  Moreover, given a modular form $f$ of level
$N$ with real coefficients, we find that the quantity
$e^{-\frac{ik\theta}{2}} f(-\frac{1}{N} + \frac{1}{N} e^{i\theta})$
is real for $0 \leq \theta \leq \pi$; this is a natural place to
look for zeros.

Our argument for the case $N=2$ relied on the ability to find a
particular horizontal contour meeting certain criteria.  If we can
find a contour passing below the points $z = -\frac{1}{N} +
\frac{1}{N} e^{i\theta}$ and $\frac{z}{Nz+1}$ and above all other
images of $z$ under the action of $\Gamma_0(N)$, then an application
of Cauchy's residue theorem to the generating function shows that
\[e^{\frac{ik\theta}{2}} e^{-\frac{2\pi n \sin \theta}{N}} f_{k, n}^{(N)}
\left(-\frac{1}{N} + \frac{1}{N} e^{i\theta}\right) - 2\cos\left(
\frac{k\theta}{2} + \frac{2\pi n}{N} - \frac{2\pi n}{N}
\cos\theta\right)\] is equal to $e^{\frac{ik\theta}{2}}
e^{-\frac{2\pi n \sin \theta}{N}}$ times the integral of the
generating function along this contour.   If we take the absolute
value of the latter expression, then choosing a contour with
horizontal height $A'$ means we can pull $e^{-\pi n (\frac{2}{N}
\sin \theta - 2A')}$ outside of the integral, and as long as
$\frac{1}{N} \sin\theta > A'$, this term has exponential decay as
$n$ grows.  This means that any bound on the size of the remaining
integral term proves a result for the location of some of the zeros
of $f_{k, n}^{(N)}(z)$ for large enough $n$.

It is possible to find such a contour for the case $N=3$.  We take a
fundamental domain in the upper half plane defined by $\{z:
-\frac{1}{2} \leq \Re(z) < \frac{1}{2} \} \cap \{z: |z-\frac{1}{3}|
> \frac{1}{3}\} \cap \{z: |z+\frac{1}{3}| \geq \frac{1}{3}\}$.  \begin{center}
\includegraphics[height=3in]{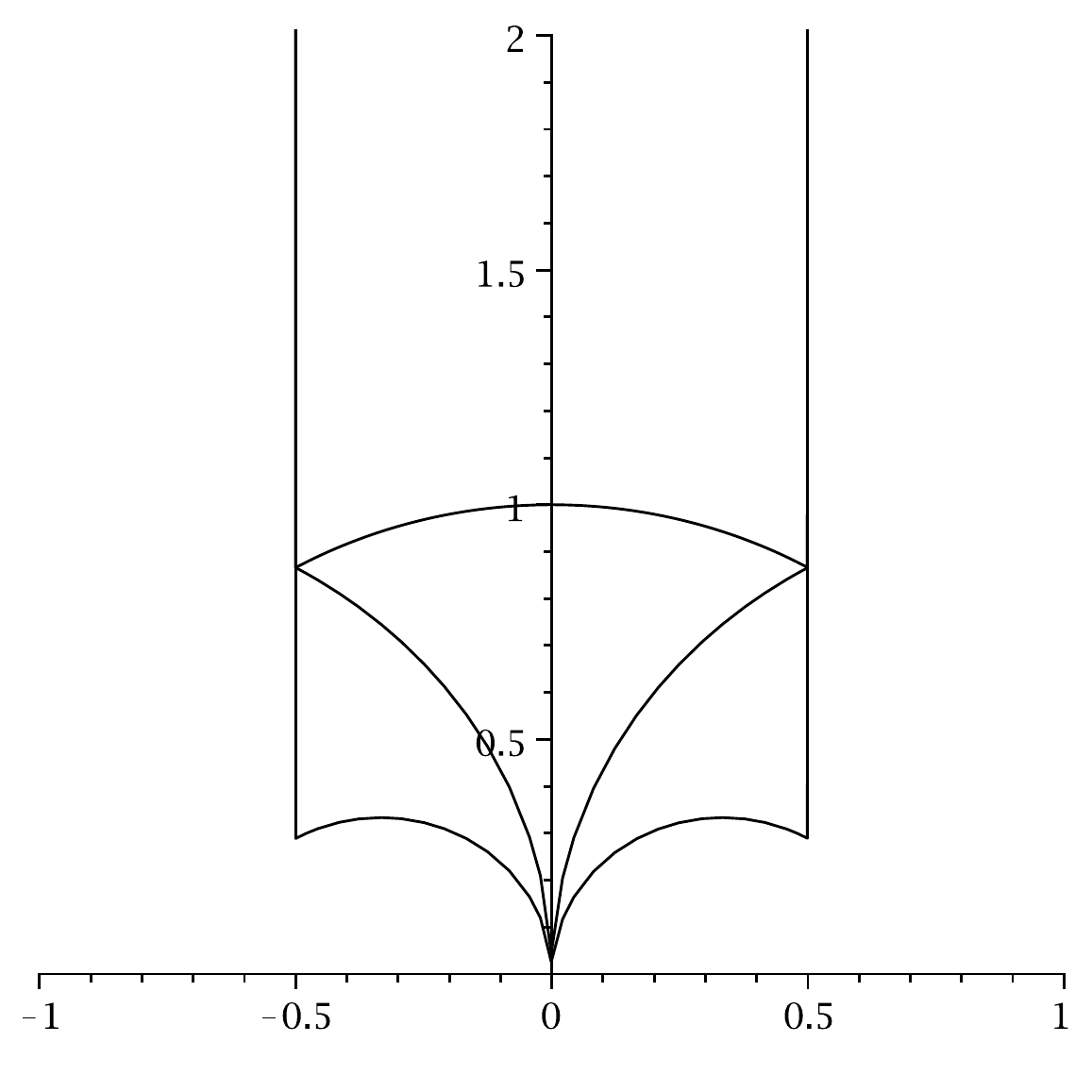}
\end{center}
We let $k = 6\ell + k'$ for $k' \in \{0, 2, 4\}$, and for $n \geq
-2\ell - \lfloor \frac{k'}{3} \rfloor$ we define a basis for
$M_k^\sharp(3)$ by letting $f_{k, n}^{(3)}(z)$ be the unique element
of $M_k^\sharp(3)$ with Fourier expansion beginning $q^{-n} +
O(q^{2\ell+\lfloor \frac{k'}{3} \rfloor + 1})$. We let $z =
-\frac{1}{3} + \frac{1}{3} e^{i\theta}$, and require the height of
the contour to be above $\frac{1}{9}$ to avoid any images of $z$
other than $\frac{z}{3z+1}$.   We restrict the range of $\theta$
values so that $z$ and $\frac{z}{3z+1}$ lie above this contour. The
contour may be as close to $\frac{1}{9}$ as we like, and we may
allow $z$ to come as close to the contour as we like; once we fix
these choices, the integral will have a finite upper bound, and the
exponential term will eventually dominate when $n \geq C\ell$ for
some constant $C$, so that the difference between the weighted
modular form $f_{k, n}^{(3)}(z)$ and the cosine function is less
than $2$.  The number of zeros of $f_{k, n}^{(3)}(z)$ in the
fundamental domain is $n + \lfloor \frac{k}{3} \rfloor$, and as the
contour approaches a height of $\frac{1}{9}$ and $z$ is allowed to
approach the contour, the number of zeros that can be proved to be
on the lower boundary approaches $.9618n + .2792k$. Thus, we have
the following theorem.
\begin{thm}
Let $f_{k, n}^{(3)}(z)$ be as above. If $n$ is large enough compared
to $\ell$, then the majority of the zeros of $f_{k, n}^{(3)}(z)$ in
the fundamental domain for $\Gamma_0(3)$ lie on the lower boundary
of the fundamental domain.
\end{thm}
A similar theorem clearly holds for the modular forms $g_{k,
n}^{(3)}(z)$.

When $N=5, 7$, or $13$, the shape of the fundamental domain is more
complicated, and choosing an appropriate contour becomes
correspondingly more difficult.  We leave this as an open problem.

\end{document}